\documentclass[final,nohypdvips,onefignum,onetabnum,smallextended]{siamart220329}

\usepackage{lipsum}
\usepackage{graphicx}
\usepackage{cleveref}
\usepackage{algorithmic}
\usepackage{svg}
\usepackage{diagbox}
\usepackage{nicefrac}
\usepackage{subcaption}
\usepackage{mathtools}
\usepackage{float}
\usepackage{amsmath}
\usepackage{amssymb}
\usepackage[10pt]{moresize}
\usepackage{graphics}
\usepackage[left=3cm, right=3cm, top=2.5cm, bottom=3cm]{geometry}
\usepackage{microtype}
\usepackage[htt]{hyphenat}
\usepackage{booktabs}

\usepackage{adjustbox}
\usepackage{multirow}
\usepackage{pifont}

\ifdefined\directlua
  \usepackage{fontspec}
\else
  \usepackage[T1]{fontenc}
  \usepackage[nomath]{lmodern}
\fi

\ifpdf
  \DeclareGraphicsExtensions{.eps,.pdf,.png,.jpg}
\else
  \DeclareGraphicsExtensions{.eps}
\fi

\newsiamremark{remark}{Remark}
\newsiamremark{hypothesis}{Hypothesis}
\crefname{hypothesis}{Hypothesis}{Hypotheses}
\newsiamthm{claim}{Claim}

\usepackage{amsopn}

\usepackage{julia-mono-listings}
\usepackage[frozencache, cachedir=.]{minted}

\usepackage[cmintegrals,cmbraces,vvarbb]{newtxmath}
\usepackage{wrapfig}
\usepackage{adpwrapfig}

\usepackage[numbers]{natbib}

\ifdefined\directlua
  \usepackage{fontspec}
\else
  \usepackage[T1]{fontenc}
  \usepackage[nomath]{lmodern}
\fi

\ifpdf
\hypersetup{
  pdftitle={},
  pdfauthor={}
}
\fi

\newcommand{\fde}{FractionalDiffEq.jl}
\newcommand{\fdesolver}{FdeSolver.jl}
\newcommand{\footnhref}[2]{\href{#1}{#2}\footnote{\url{#1}}}
\renewcommand{\descriptionlabel}[1]{\hspace*{3em}\underline{#1}}

\headers{High-Performance Fractional Differential Equation Solvers}{Q. Qu}

\title{FractionalDiffEq.\lowercase{jl}: High Performance Fractional
Differential Equation Solver in Julia\thanks{Submitted to the editors. }}

\author{Qingyu Qu\thanks{Department of Control Science and Engineering, Zhejiang University, China (\email{qingyuqu@zju.edu.cn}).},
\and
Wei Ruan\thanks{Department of Control Science and Engineering, Zhejiang University, China (\email{ruanwei@zju.edu.cn}).}
}

\begin{document}

\maketitle

\begin{abstract}
  We present \fde, a comprehensive solver suite for solving fractional differential equations, featuring high-performance numerical algorithms in the Julia programming language. \fde~is designed to be user-friendly and scalable, tackling different types of fractional differential equations, encompassing powerful numerical algorithms including predictor-corrector methods, product-integral methods, and linear multistep methods, etc, and providing a unifying API to accommodate diverse solver features. This paper illustrates the convenient usage of \fde~in modeling various scientific problems, accompanied by detailed examples and applications. \fde~leverages best practices in Julia to ensure the high performance of numerical solvers. To validate the efficiency of \fde~, we conducted extensive benchmarks that prove the superiority of \fde~against other implementations on both stiff and non-stiff problems. We further demonstrate its capability on several challenging real-life scenarios including parameter estimation in fractional-order tequila fermentation processes, and harmonic oscillator problems, etc, emphasizing the robustness and flexibility of \fde. 
\end{abstract}



\section{Introduction}\label{sec:introduction}

As modeling complexity increases in various scientific fields, it is crucial to establish an appropriate numerical model to capture the underlying dynamics of complicated systems. In this context, fractional-order models are gaining significant traction due to their superior and unique approach to depicting complex systems than integer-order derivatives. Unlike traditional methods that rely solely on first-order or second-order derivatives, fractional-order derivatives offer a special view of relationships within models by leveraging the distinctive attributes of fractional calculus. This allows for a more detailed representation of non-locality and memory effects in physical and engineering processes. Despite the evident advantages of fractional differential equations in modeling complex systems, the development of robust, general-purpose solvers for solving fractional differential equations has not kept pace with those for ordinary differential equations, and the development of related numerical solver software has lagged considerably far behind. This backwardness has hindered the advancement of numerical methods for fractional models and constrained their application in tackling complex problems across various scientific and engineering disciplines.

To fill this gap, we present \fde~- an open-source, high-performance fractional differential equation solving framework implemented in native Julia programming language \cite{bezanson2017julia}. \fde~provides a flexible, easy-to-use interface for specifying and solving fractional ordinary differential equations (FODE). We implemented efficient solvers from existing literature on fractional differential equations and conducted high-level abstraction and special optimization techniques like multiple dispatch, just-in-time (JIT) compilation, and type inferences from Julia in \fde~to accelerate the FODE numerical solving process. \fde~is designed to follow the \footnhref{https://sciml.ai/}{Scientific Machine Learning (SciML) package ecosystem of Julia} and is deeply integrated with the current differential equations solving ecosystem by C. Rackaukas et al. \cite{rackauckas2017differentialequations}. \fde~is a pioneering numerical software designed to bridge the gap between theoretical advancements in numerical methods for fractional ordinary differential equations and their practical implementation in numerical solvers for real world applications. By integrating state-of-the-art numerical techniques with high-performance computational strategies, \fde~enhances the accuracy, efficiency, and stability of fractional-order modeling, thereby expanding the applicability of fractional calculus across diverse scientific and engineering domains.

We demonstrate the capabilities of \fde~with a series of carefully designed numerical experiments. The software reliably solves meticulously chosen test problems, comparative evaluations indicate that \fde~outperforms existing alternatives, including \fdesolver, the pycaputo toolkit, and MATLAB routines, in terms of computational efficiency, accuracy, and stability, demonstrating \fde~'s robustness and superiority. Furthermore, we show \fde's applicability in real-world problems by utilizing \fde~on real-world applications like fitting fermentation fermentation dynamics with fractional differential equations in \Cref{subsec:fermentation} and modeling fractional harmonic oscillator with linear multi-term fractional differential equations in \Cref{subsec:harmonic-oscillator}. As a mature numerical software for fractional differential equations, \fde~offers both efficiency and flexibility which unlocks new modeling capabilities across various application domains. This enhancement is critical in supporting more sophisticated and accurate simulations in fields such as bioengineering, finance, and physics, where the unique capabilities of fractional models to capture memory and hereditary effects can lead to more realistic and predictive models. Thus, our solver stands as a significant advancement in the field of numerical computation, broadening the scope and depth of fractional differential equations in practical applications.

This study utilizes \fde~Version 3.5.0, it is accompanied by a permanent Zenodo DOI (https://zenodo.org/records/14313245). The accompanying user manual and documentation provide comprehensive details into the package's structure, implementation, and usage. Additionally, the user manual is also dynamically updated on the official package website (https://scifracx.github.io) to keep a continuous record of updates, new releases, and important features.

This paper is organized as follows, \Cref{subsec:existing-software} illustrates the current available toolkit for numerically solving fractional ordinary differential equations, in \Cref{sec:mathematical_description} we recall some basic definitions concerning fractional calculus and fractional ordinary differential equations solvers, and in \Cref{sec:usage_and_bencmarks}, we provide practical examples that showcase the functionality of \fde~and conducted extensive benchmarking analysis to evaluate the performance of \fde~in comparison to other prominent solvers, including \fdesolver, pycaputo, and MATLAB routines for fractional ordinary differential equations, and in \Cref{sec:applications}, we test \fde's capabilities on large fractional-order system and further prove the high-performance and usability of \fde, in \Cref{sec:discussion}, we discuss the current status in numerical solving fractional differential equations and why \fde~can be helpful for numerical modeling and simulation of fractional-order systems, in \Cref{sec:conclusion} we summarize the key findings of this study, conclude what this study can bring in for the research of fractional differential equations, and outline the further development directions for \fde~.

\subsection{Comparison to Existing Software}\label{subsec:existing-software}

Numerical solvers for fractional differential equations face significant challenges, particularly in terms of flexibility and computational efficiency. Compared to ordinary differential equation solvers, they lack the advanced features and optimized performance. Existing algorithms for fractional ordinary differential equations are often plagued by inefficiencies and lack the versatility required to effectively address a wide variety of complex problems. Consequently, the development of numerical solvers for fractional differential equations lags significantly behind comparing with ODE numerical solvers. \cite{li2017review}.

R. Garrappa et al. \cite{garrappa2018numerical, garrappa2020initial} developed a suite of MATLAB routines for numerically solving fractional ordinary differential equations employing product integral rules and fractional linear multistep methods, including both system of fractional differential equations and linear multi-term fractional ordinary differential equations. Numerical algorithms in \fde~for numerically solving fractional ordinary differential equations are primarily inspired by the methodologies proposed by R. Garrappa \cite{garrappa2018numerical}. These algorithms have been further refined and tailored to leverage the specific features of the Julia programming language, incorporating optimizations that enhance computational efficiency, numerical stability, and overall usability. By exploiting Julia’s just-in-time (JIT) compilation, multiple dispatch, and efficient memory management, these implementations achieve significant performance improvements while maintaining high accuracy. Additionally, specialized data structures and parallel computing techniques are employed to further accelerate computations, making \fde~a robust and user-friendly framework for researchers and practitioners working with FODEs.

I. Podlubny et al. \cite{podlubny2000matrix} developed a matrix discretization method in which fractional-order derivatives and fractional-order integrals are represented by an upper and lower triangular strip matrix with special coefficients in fractional differential equations. By substituting the differential operators with such matrices, the discretized fractional differential equations are thus transformed into a large algebraic equation where computing the final numerical solution is just to solve the large linear system. Though the Matrix-Discretizition methods make solving fractional differential equations intuitive and intelligible, the Matrix-Discretizition methods assume the zeros initial conditions of FODE, and the nonzero initial value FODE needs to be transformed into zero conditions initial value FODE. Furthermore, the benchmarking presented in \Cref{fig:linear-multiterms-benchmarks} indicates that the Matrix-Discretizition methods do not exhibit superior performance compared to other iterative numerical solvers. In particular, while matrix-based approaches offer a structured way to approximate numerical solutions, they often suffer from higher computational complexity and memory requirements when dealing with fine-grained mesh for large-scale FODE system.

M. Khalighi et al. \cite{fdesolver} developed a similar toolbox \fdesolver~for solving fractional ordinary differential equations in Julia problem-solving environment, in which they implemented benchmarks between \fdesolver~and \fde~and MATLAB routines to compare the superiority of \fdesolver. In \fdesolver, there are only the Product Integral method and Newton-Raphson method for fractional ordinary differential equations, which don't have mature implementations for numerical solvers of linear multi-term fractional differential equations and their implementation doesn't have a judging criterion that indicates the successful solving process, which makes it troublesome to assess whether the previous iteration has a successful convergence. Besides, when handling stiff fractional ordinary differential equations, their solvers suffer divergence and fail to provide accurate results. In the benchmarks of \cite{fdesolver} which compares \fde, \fdesolver~and MATLAB routines, the benchmarks are biased and don’t reflect the real performance of \fde, and the related paper is not rigorous in specifying their testing problems. In addition, benchmarks in \cite{fdesolver} is only limited to three sources of FODE solvers, doesn't include the implicit methods for accurate numerical solutions. This paper will not only represent the superior numerical in \fde~, but also discuss a detailed and thorough comparison between these numerical solvers and provide a fair comparison of the current existing FDE solvers on standard test problems.

A. Fikl et al. \cite{alex_fikl_2024_13926507} implemented a Python software pycaputo for fractional calculus evaluation and fractional differential equations solving. Pycaputo supports fractional derivatives and integrals evaluating and incorporates predictor-corrector methods, trapezoidal methods, and Euler methods for numerically solving fractional ordinary differential equations in the Caputo sense. It should be noted that \cite{alex_fikl_2024_13926507} has an adaptive time stepping strategy which is based on \cite{jannelli2020novel}. However, the problem types that pycaputo supports are limited and its numerical solvers' performance is relatively bad compared to its Julia and MATLAB equivalents.

D. Xue et al. developed the FOTF toolbox \cite{xue2019fotf, xue2024fractional} with a similar interface of MATLAB control system package for fractional-order control system modeling. While FOTF is primarily designed for fractional order control system modeling, it is noteworthy that there are also system of fractional ordinary differential equations and linear multi-term fractional ordinary differential equations numerical solvers implemented in the FOTF toolbox. Specifically, the toolbox includes \texttt{fode\_caputo9} for non-zeros initial conditions linear multi-term fractional ordinary differential equations which employ the Grunwald-Letnikov difference scheme and \texttt{nlfode\_vec} for nonlinear fractional ordinary differential equations systems which also utilized Grunwald-Letnikov difference schemes. In this paper, we only compare the numerical solvers \texttt{nlfode\_vec} which is capable of handling non-zero initial conditions FODE problems and \texttt{fode\_caputo9} for non-zero initial conditions of linear multi-term FODEs.

In contrast to existing numerical software, \fde~is designed to be robust, performant and scalable while having controls to be customized to specific applications. \fde~has the most extensive numerical solvers available for tackling different kinds of fractional differential equations. Each solver is designed to be computationally efficient and numerically accurate to support complex and high-demanding modeling tasks. Convenient special functions handling is important when solving fractional differential equations, \fde~integrate fast Mittag-Leffler function computing from the algorithm of \cite{garrappa2015numerical} and provide a user-friendly interface to facilitate users with hands-on Mittag-Leffler function computing.

\section{Mathematical Description}\label{sec:mathematical_description}

This section introduces the primary mathematical definitions of fractional calculus and fractional differential equations and recalls the properties that will be used in the subsequent sections. For a more detailed and comprehensive introduction to fraction calculus and fractional differential equations, we recommend readers to refer to textbooks \cite{xue2024fractional,podlubny1998fractional,oldham1974fractional,singh2022handbook,baleanu2012fractional,gorenflo1997fractional,li2015numerical} and related papers \cite{garrappa2018numerical,blank1996numerical}.

Suppose $D^\alpha_t$ is Caputo sense time fractional derivative with order $\alpha\in[0,1]$. The Caputo sense fractional-order derivative in a sufficient differentiable function is given by:
\begin{equation}
    {^C_{t_0}D^\alpha_t}f(t) = I_{t_0}^{m-\alpha}f^{(m)}(t) = \frac{1}{\Gamma(m-\alpha)}\int^t_0\frac{f^{(m)}(\psi)d\psi}{(x-\psi)^{\alpha-m+1}},\ m-1<\alpha<m
\end{equation}
where $I_{t_0}^{m-\alpha}$ is Riemann-Liouville sense time fractional integral of order $m-\alpha$ which is defined as:
\begin{equation}
    I^{\beta}_{t_0}g(t)=\frac{1}{\Gamma(\beta)}\int_{t_0}^t\frac{g(\psi)d\psi}{(t-\psi)^{1-\beta}}
\end{equation}
where $f,g : \mathbb{R}^n \times \mathbb{P} \to \mathbb{R}^n$, $m = \lceil\alpha\rceil$ is the smallest integer greater or equal to $\alpha$. In this paper, we will assume that $f$ is continuous and differentiable. With Caputo sense fractional-order derivative, an initial value problem for a fractional-order ordinary differential equation system in Caputo sense can be formulated as :
\begin{equation}
    \begin{cases}
        {^C_0D_t^\alpha y(t) = f(t, y(t))}\\
        y(t_0) = y_0,\ y'(t_0)=y_0^{(1)},\cdots,y^{(m-1)(t_0)}=y_0^{(m-1)}
    \end{cases}
    \label{eq:fractional-differential-equations}
\end{equation}
where the right-hand side function $f(t,y(t))$ is assumed to be continuous in the time span $t\in[t_0, t_f]$.

A special case of linear fractional-order ordinary differential equations often occurs in fractional-order modeling, especially when one encounters transfer functions in fractional-order control system \cite{xue2024fractional}, and in this case, the linear fractional ordinary differential equations can be formulated as linear multi-term fractional-order ordinary differential equations, which have the formulation of:
\begin{equation}
    \lambda_Q{^CD_{t_0}^{\alpha_Q}}y(t) + \lambda_{Q-1}{^CD_{t_0}^{\alpha_{Q-1}}}y(t) + \cdots + \lambda_2{^CD_{t_0}^{\alpha_2}}y(t) + \lambda_1{^CD_{t_0}^{\alpha_1}}y(t) = f(t, y(y))
\end{equation}
where the coefficients $\lambda_1,\cdots\lambda_Q$ and orders $\alpha_1,\cdots,\alpha_Q$ are assumed to be sorted in ascending order, but $f(t, y(t))$ could be a nonlinear function of $y(t)$. For a more rigorous and detailed discussion on the topic of fractional calculus and fractional differential equations, we refer the readers to \cite{podlubny1998fractional,oldham1974fractional,li2015numerical,doi:10.1142/10044}.

\subsection{Fractional Linear Multiple Step Methods}

Compared to one-step methods, multistep methods incorporate previous approximations to achieve higher precision solutions. And because of the special memory effect(i.e. non-locality in time) of fractional derivatives, unlike the fixed number of steps used in multistep methods in ODE, multistep methods in FODE would not have fixed steps, and the whole memory needed to be considered. Fractional linear multistep methods(FLMM) have been extensively studied by \cite{lubich1986discretized}, the core idea of FLMM is to utilize the convolution quadrature for the Riemann-Liouville integral.
\begin{equation}
    y(t_n)=\sum_{k=0}^{m-1}\frac{t_n^k}{k!}y_0^{(k)} + h^\alpha\sum_{i=0}^s \omega_{n,i}f(t_i, y(t_i)) + h^\alpha\sum_{i=0}^n\omega_{n-i}^{(\alpha)}f(t_i, y(t_i))
\end{equation}
on uniform grids $t_n=t_0+nh$. While $\omega_n^{(\alpha)}$ are depend on the characteristic polynomials $\rho(z)$ and $\sigma(z)$, it can be calculated more easily using the convolution quadrature weights obtained from:
\begin{equation}
\displaystyle\sum_{n=0}^\infty\omega_n^{(\alpha)}\psi^n = \omega^{(\alpha)}(\psi),\quad \omega^{(\alpha)}(\psi) = (\delta(\psi))^{-\alpha}
\end{equation}
where the generating function for the fractional linear multistep method is derived from G. Roberto \cite{garrappa2015trapezoidal}, which provides a comprehensive framework for constructing numerical methods for fractional differential equations based on convolution quadrature. \Cref{tab:flmm_methods} represents the necessary generating functions for FLMM methods:

\begin{table}[H]
    \centering
    \begin{tabular}{cl}
    \toprule
         \textbf{FLMM method} & \textbf{Generating function} \\
         \midrule
         BDF2 & $\delta(\xi) = \frac{\sigma(1/\xi)}{\rho(1/\xi} = \frac{3}{2}-2\xi+\frac{1}{2}\xi^2$ \\
         Trapezoidal & $\delta(\xi) = \frac{\sigma(1/\xi)}{\rho(1/\xi} = \frac{1}{2}(1+2\displaystyle\sum_{n=1}^\infty\xi^n)$ \\
         Newton Gregory & $\delta_\alpha(\xi) = \frac{\sigma(1/\xi)}{\rho(1/\xi} = (1-\xi)^{-\alpha}(1-\frac{\alpha}{2}(1-\xi))$ \\
         \bottomrule
    \end{tabular}
    \caption{A list of FLMM methods and their generating functions $\delta(\xi)$.}
    \label{tab:flmm_methods}
\end{table}

\subsection{Product Integration rules}

Fractional ordinary differential equations in \Cref{eq:fractional-differential-equations} can be reformulated as weakly-singular Volterra integral equations by applying a Riemann-Liouville integrator to both sides of the equation. This transformation allows for a reformulation that facilitates numerical treatment of fractional derivatives. To discretize these equations, we utilized the product integration rules introduced by Young \cite{young1954approximate}, which are specifically designed to handle second-kind weakly-singular VIEs on the discretized mesh. By applying these rules, the discretized fractional ordinary differential equations can be reformulated in a piecewise way \cite{garrappa2018numerical}:
\begin{equation}
    y(t_n)=\sum_{k=0}^{m-1}\frac{t_n^k}{k!}y_0^{(k)}+\frac{1}{\Gamma(\alpha)}\sum^{n-1}_{j=0}\int^{t_{j+1}}_{t_j}(t_n-\tau)^{\alpha-1}f(\tau, y(\tau))d\tau
    \label{eq:weakly-singular-vie}
\end{equation}
then based on different numerical approximation strategies for the integral term in \Cref{eq:weakly-singular-vie}, a variety of explicit and implicit FODE integration methods can be derived. These methods offer different trade-offs in terms of accuracy, stability, and computational efficiency, providing flexibility in solving fractional-order systems across various applications.

\begin{equation}
    \text{Implicit Trapezoidal}: y(t_n)=\sum_{k=0}^{m-1}\frac{t_n^k}{k!}y_0^{(k)}+\frac{1}{\Gamma(\alpha+2)}(\tilde{a}_{\alpha,n} + \sum^{n}_{i=1}a_{\alpha,n}f(t_i, y_i))
    \label{eq:implicit-pi-trapezoidal}
\end{equation}

\begin{equation}
    \text{Implicit Rectangular}: y(t_n)=\sum_{k=0}^{m-1}\frac{t_n^k}{k!}y_0^{(k)}+\frac{1}{\Gamma(\alpha+1)}\sum^{n-1}_{i=0}b_{\alpha,n-i-1}f(t_i, y_i)
    \label{eq:implicit-pi-rectangular}
\end{equation}

\begin{equation}
    \text{Explicit Trapezoidal}: y(t_n)=\sum_{k=0}^{m-1}\frac{t_n^k}{k!}y_0^{(k)}+\frac{1}{\Gamma(\alpha+1)}\sum^{n}_{i=0}b_{\alpha,n-i}f(t_i, y_i)
    \label{eq:explicit-pi-rectangular}
\end{equation}
with the computing coefficients:

\begin{equation}
    b_{\alpha,n}=(n+1)^\alpha-n^\alpha,\ \tilde{a}_{\alpha,n}=(n-1)^{\alpha+1}-n^\alpha(n-\alpha-1)
\end{equation}

\begin{equation}
    a_{\alpha,n}=(n-1)^{\alpha+1}-n^\alpha(n-\alpha-1)-2n^{\alpha+1}+(n+1)^{\alpha+1}
\end{equation}

\subsection{Predictor-Corrector Methods}

Predictor-corrector methods for fractional ordinary differential equations were first proposed by K. Diethelm et al. \cite{diethelm2002predictor} and have since gotten significant attention. As a generalization of the Adams-Bashforth-Moulton multistep method, the predictor-corrector method starts with an estimation of the unknown value as the \texttt{predictor}, then the \texttt{corrector} would refine it in specified steps. The core idea behind the predictor-corrector method is to first generate an initial estimate of the unknown function value using a predictor step, followed by a corrective refinement through the corrector step. So algorithms can be simplified as a predictor:
\vspace{-2em}
    \begin{align}
        & y^P(t_n) = \sum_{k=0}^{m-1}\frac{t_n^k}{k!}y_0^{(k)} + h^\alpha\sum_{i=0}^{n-1}b_{n-i-1}^\alpha f(t_i, y(t_i))
    \end{align}
    \vspace{-2em}
and a corrector:
\vspace{-2em}
\begin{align}
    & y(t_n) = \sum_{k=0}^{m-1}\frac{t_n^k}{k!}y_0^{(k)} + \frac{h^\alpha}{\Gamma(\alpha+2)}f(t_n, y^P(t_n)) + \frac{h^\alpha}{\Gamma(\alpha+2)}\sum_{i=0}^{n-1}a_{n-i-1}f(t_i, y(t_i))
\end{align}
while there convolution sum in both corrector and predictor, and such calculation could exploit fast Fourier transform which can reduce the computational cost to $\mathcal{O}(n(\log_2n)^2)$.

\section{Usage and Benchmarks}\label{sec:usage_and_bencmarks}

In the following section, we provide a comprehensive guide on utilizing \fde~through a series of numerical examples. These examples are selected to illustrate the tool’s intuitive usability, highlighting its flexibility and computational robustness in diverse scenarios. By systematically exploring various problem settings, we aim to highlight the solver's versatility in addressing various kinds fractional differential equations. Furthermore, we conducted exhaustive benchmarking analyses to compare the performance of \fde~with the solvers discussed in the previous section. These benchmarks evaluate key metrics, including computation speed and accuracy, across a range of test cases that contain both non-stiff and stiff problems. By systematically presenting these results, we aim to provide a clear and objective assessment of \fde’s capabilities. The findings from these benchmarks underscore the superior performance of \fde~compared to other solvers. Its ability to handle complex fractional differential problems with exceptional speed and precision not only validates its computational advantage but also establishes its position as a cutting-edge tool in the field. This section thus serves a dual purpose: demonstrating the practical application of \fde~while providing robust evidence of its superiority over existing solutions.

\subsection{Intuitive \fde~usage}

In the following subsection, we delve into the practical application of \fde~by introducing a series of classical numerical examples involving system fractional ordinary differential equations (FODEs) and linear multi-term fractional ordinary differential equations. These examples are carefully chosen to represent a wide range of scenarios, highlighting the versatility and robustness of \fde~in solving complex problems in fractional calculus. The basic structure of \fde~when solving fractional differential equations with different types of problem is shown in \Cref{fig:architecture}:

\begin{figure}[H]
    \centering
    \includegraphics[width=\textwidth]{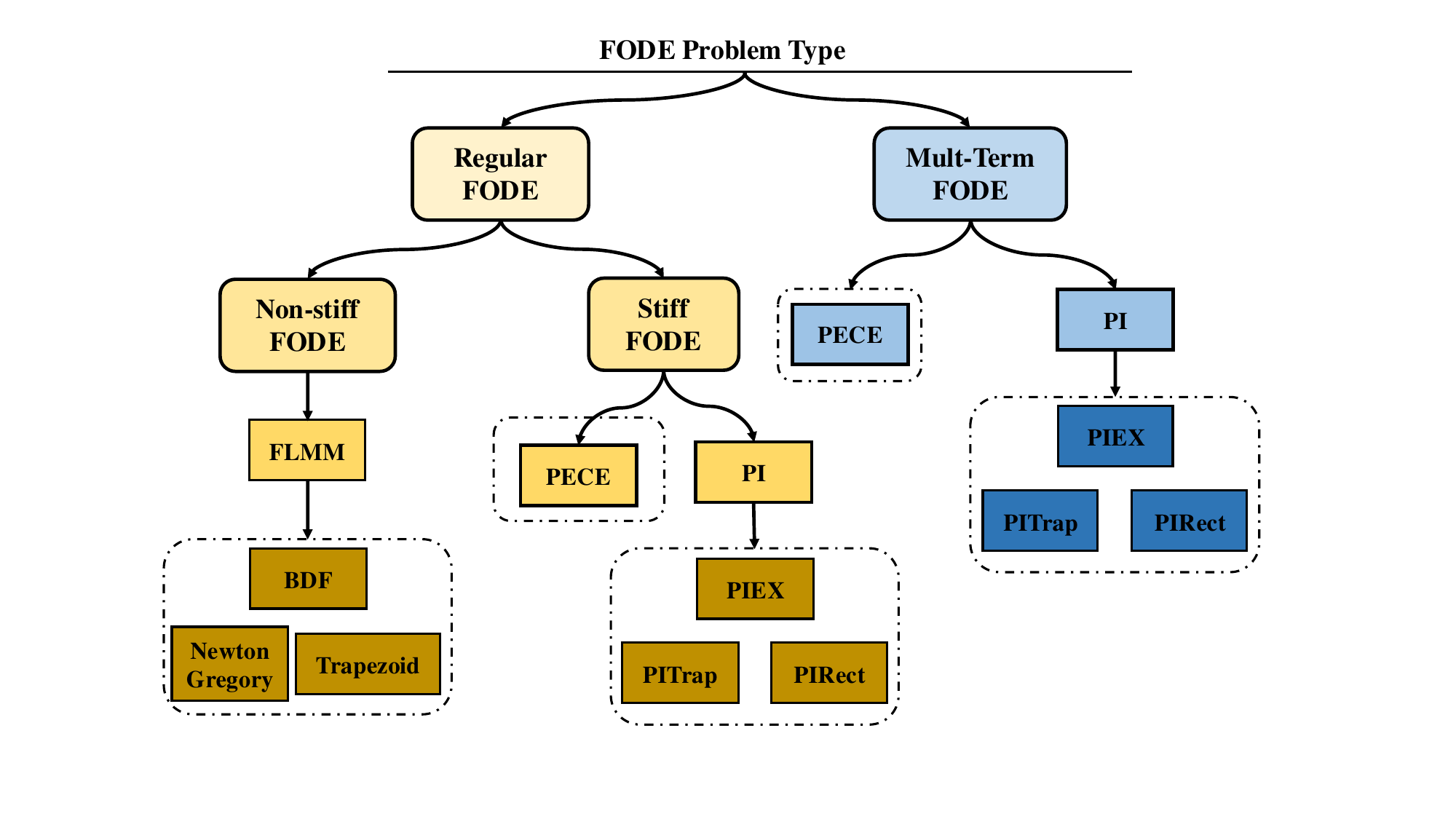}
    \caption{\textbf{Schematic of \fde~:} \fde~has a modular architecture where different algorithms should be chosen concerning the problem properties according to different kinds of fractional ordinary differential equations.}\label{fig:architecture}
\end{figure}

\subsubsection{Solve fractional order Chua chaotic system}\label{subsubsec:solve_fode_chaotic_system}

Classical Chua's circuit is a simple model that exhibits nonlinear dynamical phenomena including bifurcation and chaos with appropriate configuration. According to I. Petras et al. \cite{petravs2011fractional}, when the capacitor and inductor in Chua's system are replaced by the new linear capacitor model and inductor model which is based on Curie's empirical law of 1889 \cite{westerlund}, the general input voltage and current would behave a fractional-order derivative relationship as:
\begin{equation}
    I(t)=C\frac{d^\alpha V(t)}{dt^\alpha}\equiv C_0D^\alpha_tV(t),\ V(t)=L\frac{d^\alpha I(t)}{dt^\alpha}\equiv L_0D^\alpha_tI(t)
\end{equation}
where $C$ is the capacitance of the real capacitor and $L$ is the inductance of the inductor. The Caputo derivative here can better describe the memory effects of the charging and discharging process of capacitor, and the energizing and de-energizing process of inductor \cite{westerlund1991dead}. Hence we have the fractional-order Chua chaotic system as:
\begin{equation}
    \begin{cases}
        {^C_0D^{\alpha_1}_t} u_1(t)= a(u_2(t)-u_1(t)-f(u_1))\\
    {^C_0D^{\alpha_2}_t} u_2(t) = u_1(t)-u_2(t)+u_3(t)\\
    {^C_0D^{\alpha_3}_t} u_3(t) = -bu_2(t)-cu_3(t)\\
    \end{cases}
\end{equation}

To utilize \fde~solve the fractional Chua system, we simply need to follow the use the "define the problem then solve the problem" convention. To maintain consistency with the interface design of \cite{rackauckas2017differentialequations} and ensure seamless integration within the SciML differential equations solvers ecosystem, \fde~adopts a similar formulation style for defining fractional ordinary differential equations. This design choice facilitates user familiarity and interoperability with existing numerical solvers, thereby enhancing usability and accessibility for researchers and practitioners. The primary distinction, however, lies in the explicit specification of the fractional order in the Caputo sense, which is a required parameter in \fde~during problem construction.

\begin{figure}[H]
    \hspace*{0.3in}
    \begin{minipage}[c]{0.4\textwidth}
        \centering
        \begin{minted}[breaklines,escapeinside=||,mathescape=true, linenos, numbersep=3pt, gobble=2, frame=lines,fontsize=\small, framesep=2mm]{julia}
using FractionalDiffEq
function chua!(du, x, p, t)
    a, b, c, m0, m1 = p
    du[1] = a*(x[2]-x[1]-(m1*x[1]+m0*(abs(x[1]+1)-abs(x[1]-1))))
    du[2] = x[1]-x[2]+x[3]
    du[3] = -b*x[2]-c*x[3]
end
alpha = [0.93, 0.99, 0.92];
u0 = [0.2; -0.1; 0.1];
tspan = (0, 100);
p = [10.725, 10.593, 0.268, -0.1927, -0.7872]
prob = FODEProblem(chua!, alpha, u0, tspan, p)
sol = solve(prob, BDF(), dt=0.01)
        \end{minted}
\end{minipage}
    \hfill
    \begin{minipage}[c]{.50\textwidth}
        \centering
        \includegraphics[width=0.9\textwidth]{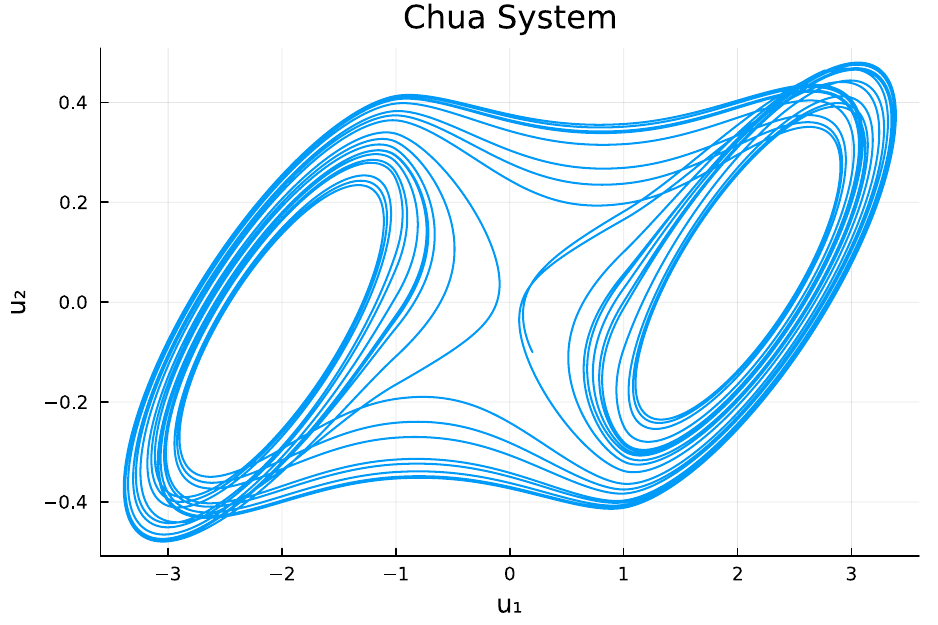}
    \end{minipage}
    \caption{\fde~allows seamless solvers switching between different FODE numerical algorithms for FODE problems. We can choose different solvers according to the commensurate order of the FODE system from available solvers in \Cref{tab:fode_solvers}.}\label{fig:fode-usage}
    \vspace{-2em}
\end{figure}

\subsubsection{Solve Bagley-Torvik equation}\label{subsubsec:solve_bagley_torvik_equation}

The Bagley-Torvik equation arises in the modeling of the motion of a rigid plate immersed in a Newtonian fluid \cite{podlubny1998fractional,diethelm2002numerical}, which describes the relationship between the stress and velocity for the geometry and loading can be described in a fractional-order manner:
\begin{equation}
    \sigma(t,z)=\sqrt{\mu\rho}{^C_0D^{1/2}_t}v(t,z)
\end{equation}
where $\rho$ is the fluid density, $\mu$ is the viscosity of the fluid and . So the dominant fractional differential equation can be formulated:
\begin{equation}
\begin{aligned}
    Au''(t)+&B{^C_0D^{3/2}_t}u(t)+Cu(t)=f(t)\\
    &f(t) = \begin{cases}
        0,\ t>1\\
        8,\ 0\leq t\leq 1
    \end{cases}
    \end{aligned}
    \label{eq:bagley-torvik-equation}
\end{equation}
where $A\neq0$ and $B,C \in \mathcal{R}$, the initial conditions are chosen as $u(0)=u'(0)=1$ since the displacement and velocity of the plate-fluid system must be in equilibrium for the fractional derivative to be applicable \cite{podlubny1998fractional}. The fractional derivative $D^{3/2}_tu(t)$ introduced in this equation is used to interpret the relationship of stress at a given point at any time is dependent on the time history of the velocity profile at the point.

\begin{figure}[H]
\hspace*{0.3in}
    \begin{minipage}[c]{0.4\textwidth}
        \centering
        \begin{minted}[breaklines,escapeinside=||,mathescape=true, linenos, numbersep=3pt, gobble=2, frame=lines, fontsize=\small, framesep=2mm]{julia}
using FractionalDiffEq

tspan = (0.0, 20.0)
rhs(u, p, t) = ifelse(t > 1, 0, 8)
u0 = [0.0, 0.0]
prob = MultiTermsFODEProblem([1, 1/2, 1/2], [2, 1.5, 0], rhs, u0, tspan)
sol = solve(prob, MTPECE(), dt=0.01)
        \end{minted}
\end{minipage}
    \hfill
    \begin{minipage}[c]{.50\textwidth}
        \centering
        \includegraphics[width=0.9\textwidth]{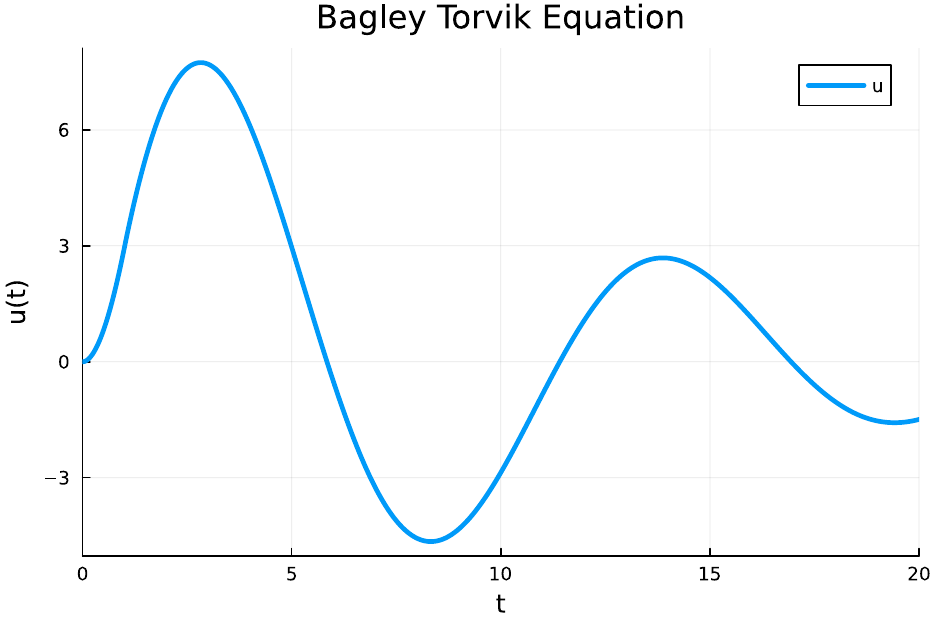}
    \end{minipage}
    \caption{\fde~allows seamless solvers switching between different numerical algorithms for multi-term FODE problems. We can choose any available multi-term FODE solver during calling \texttt{sovle} from available solvers in \Cref{tab:mtfode_solvers}.}\label{fig:multiterm-fode-usage}
\end{figure}

\subsection{Benchmarks with other FDE solvers}\label{subsec:benchmarks}

In the following subsection, we present the performance comparison between different available fractional ordinary differential equation solvers across Julia, MATLAB, and Python ecosystems.This analysis aims to provide an exhaustive and fair evaluation by systematically solving a series of well-established benchmarking problems proposed by \cite{garrappa2018numerical, XueBai+2017+1305+1312}, and \cite{implicit}.

The detailed benchmarking scripts have been made publicly available on GitHub, incorporating well documented, comprehensive and rigorously standardized experimental setups to ensure transparency and reproducibility\footnote{\href{https://github.com/SciFracX/fde_benchmarks}{https://github.com/SciFracX/fde\_benchmarks}}. All the experiments were carried out in Julia Version 1.11.0, MATLAB Version 9.12.0.1975300 (R2024a) Update 3, and Python Version 3.13.0 on a computer equipped with a CPU Intel i7-9750H at 2.60 GHz running on the MacOS 14.6 Sonoma M3 chip. To maintain integrity and fairness in our evaluation, we performed multiple runs of each benchmark using BenchmarkTools.jl \cite{revels2017benchmarktools} and provided a comprehensive user manual detailing the benchmarking procedures, facilitating the reproducibility of the evaluations. It is important to note that variations in software versions or hardware configurations may lead to discrepancies in performance outcomes.

\subsubsection{FODE solvers benchmarks}\label{subsubsec:fode_benchmarks}

To provide a comprehensive understanding of the numerical solvers evaluated in the fractional ordinary differential equation (FODE) benchmarks, we first present an overview of all participating solvers. The complete list of common fractional ordinary differential equation solvers, along with their respective sources, is systematically summarized in \Cref{tab:fode_solvers}. 
\begin{table}[H]
    \centering
    \begin{tabular}{ccccc}
    \toprule
         \fde~ & MATLAB & FdeSolver.jl & pycaputo \\
         \midrule
         PECE & FDE\_PI12\_PC & PC & PECE \\
         PIRect & FDE\_PI1\_IM & - & -  \\
         PITrap & FDE\_PI2\_IM & - & - \\
         PIEX & FDE\_PI1\_EX & - & - \\
         BDF& FLMM2(BDF) & - & -  \\
         NewtonGregory & FLMM2(NewtonGregory) & - & -  \\
         Trapezoid & FLMM2(Trapzoid) & - & -  \\
         - & FOTF(nlfode\_vec) & - & -  \\
         \bottomrule
    \end{tabular}
    \caption{A list of numerical solvers for fractional ordinary differential equations, solvers on the same row are built on the same numerical algorithm.}
    \label{tab:fode_solvers}
\end{table}
And all the available numerical solvers for linear multi-term FODE can be summarized in solvers overview in \Cref{tab:mtfode_solvers} along with their sources:
\begin{table}[H]
    \centering
    \begin{tabular}{ccccc}
    \toprule
         \fde~ & MATLAB \\
         \midrule
         MTPIEX & MT\_FDE\_PI1\_Ex   \\
         MTPIRect & MT\_FDE\_PI1\_Im  \\
         MTPITrap & MT\_FDE\_PI2\_Im   \\
         MTPECE & MT\_FDE\_PI12\_PC  \\
         - & fode\_caputo9  \\
         - & Matrix-Discretization \\
         \bottomrule
    \end{tabular}
    \caption{A list of numerical solvers for linear multi-term fractional ordinary differential equations, solvers on the same row are built on the same numerical algorithm.}
    \label{tab:mtfode_solvers}
\end{table}
Firstly, we conducted a comprehensive benchmarking analysis of all available FODE solvers by evaluating their performance on simple single-term linear and nonlinear FODEs, selected from the established test problems of \cite{garrappa2018numerical} and \cite{fdesolver}.
\begin{equation}
\begin{aligned}
    &\textbf{Linear FODE problem:} \quad 
    \mathcal{D}_{t_0}^{\alpha} u(t) =
    \frac{40320}{\Gamma(9-\alpha)}t^{8-\alpha}
    -3 \frac{\Gamma(5+\beta/2)}{\Gamma(5-\alpha/2)}t^{4-\alpha/2} 
    +\frac{9}{4}\Gamma(\alpha+1) \\
    &\quad +\left( \frac{3}{2}t^{\alpha/2}-t^{4}\right)^3 
    -\left(u(t)\right)^{3/2}, \quad u(0) = 0. \\[6pt]
    &\textbf{Analytical solution:} \quad 
    u(t) = t^8 -3t^{4+\alpha/2} +\frac{9}{4}t^\alpha.
\end{aligned}
\label{eq:single-term-linear-fode-equation}
\end{equation}

\begin{equation}
\begin{aligned}
    &\textbf{Nonlinear FODE problem:} \quad 
    \mathcal{D}_{t_0}^{\alpha}u(t) = -10u(t), \quad u(0) = 1. \\[6pt]
    &\textbf{Analytical solution:} \quad 
    u(t) = E_\alpha(-10t^\alpha).
\end{aligned}
\label{eq:single-term-nonlinear-fode-equation}
\end{equation}
where $E_{\alpha, \beta}^\gamma(z)=\frac{1}{\Gamma(\gamma)}\displaystyle\sum^\infty_{k=0}\frac{\Gamma(\gamma+k)z^k}{k!\Gamma(\alpha k+\beta)}$ is the generalized three-parametric Mittag-Leffler function, which has the relationship of $E_{\alpha}(z)=E_{\alpha, 1}(z)$, and $E_{\alpha, \beta}(z) = E_{\alpha, \beta}^1(z)$ \cite{prabhakar1971singular}. \fde~has built-in numerical Mittag-Leffler function computing which can be called using \texttt{mittleff(alpha, beta, z)} for the two-parametric Mittag-Leffler function, or \texttt{mittleff(alpha, beta, gamma, z)} for the three-parametric Mittag-Leffler function \cite{garrappa2015numerical}.

\begin{figure}[H]
    \centering
    \adjustbox{trim=4em 2em 3em 3em, clip}{
        \includegraphics[width=1.15\textwidth]{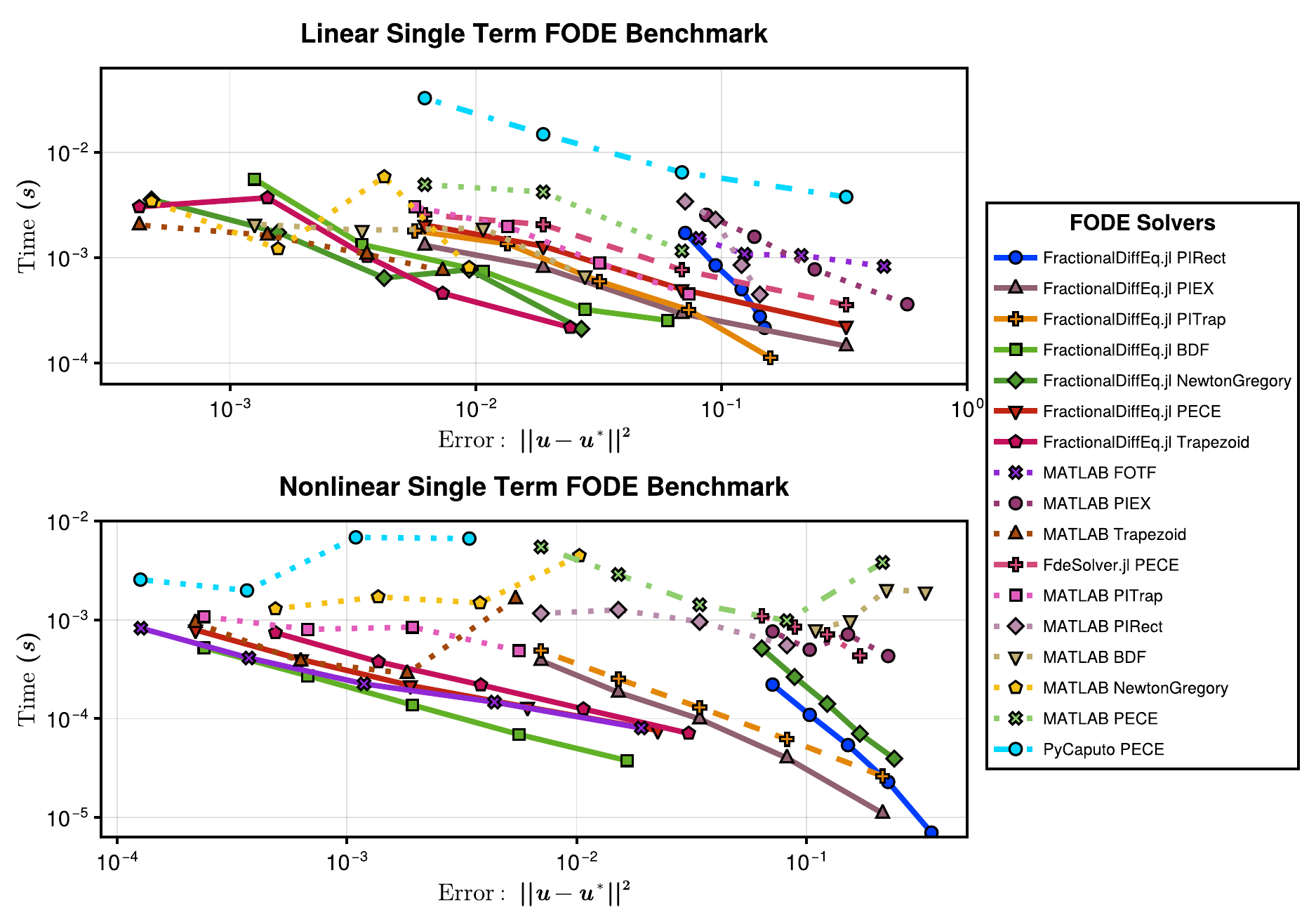}
    }
    \caption{\textbf{Work Precision Diagram of linear and nonlinear single term FODE problem}: In the benchmarking analysis of linear problema presented \cite{garrappa2018numerical}, all FODE solvers demonstrated convergence within a reasonable computational time. However, solvers from \fde~exhibited significantly performance, consistently achieving faster and more reliable solutions for both linear and nonlinear fractional ordinary differential equations compared to the corresponding MATLAB, Julia and Python counterparts.}\label{fig:singleterm-benchmarks}
\end{figure}

In \Cref{fig:singleterm-benchmarks}, we separately compare the computational efficiency of all available FODE solvers when applied to linear \Cref{eq:single-term-linear-fode-equation} and nonlinear FODE problem \Cref{eq:single-term-nonlinear-fode-equation}. The benchmarking problem set are solved with step size of $h=2^{-n}, n=2,\cdots,7$. \Cref{fig:singleterm-benchmarks} shows the efficiency comparison between different solver framework on low dimensional FODE problems. Our experiment results indicate that for relatively small and low-dimensional fractional ordinary differential equations, \fde~and MATLAB achieved an excellent balance of computational speed and numerical precision In contrast, while solvers from \fdesolver~and pycaputo toolkit being able to successfully solve the problem, they fail to efficiently handle the aforementioned FODE problems. The observed inefficiencies in these solvers suggest potential limitations in their numerical stability, or underlying discretization methods, which hinder their performance when dealing with more complex fractional differential equations.

Secondly, we separately benchmarked FODE solvers from \fdesolver, \fde, routines from MATLAB and pycaputo on benchmark problem 5 from \cite{XueBai+2017+1305+1312}, which is a non-stiff and non-commensurate nonlinear fractional-order state-space model in Caputo sense fractional derivative, the non-stiff problem in time span $t\in[0, 5]$ with initial conditions $u_0=[1, 0.5, 0.3]^\top$ can be formulated as:
\begin{equation}
    \begin{cases}
    {^C_0D}^{0.5}_t u_1(t) = \frac{1}{\sqrt{\pi}}(\sqrt[6]{(u_2(t)-0.5)(u_3(t)-0.3)} + \sqrt{t}) \\
    {^C_0D}^{0.2}_t u_2(t) = \Gamma(2.2)(u_1(t)-1)\\
    {^C_0D}^{0.6}_t u_3(t) = \frac{\Gamma(2.8)}{\Gamma(2.2)}(u_2(t)-0.5)
\end{cases}
\label{eq:nonstiff-fode-equation}
\end{equation}
where the non-stiff FODE problem has an analytical solution of $u(t) = [t+1, t^{1.2}+0.5, t^{1.8}+0.3]^\top$. The chosen non-stiff FODE problem has a solution profile where the third state $u_3(t)$ increase with $t$ rapidly, 

The stiff nonlinear FODE problems are chosen from \cite{implicit} which is a commensurate order fractional ordinary differential equation problem and can be formulated as:
\begin{equation}
    {^C_0\textbf{D}}^{\alpha}_t \textbf{u}(t) = (\textbf{A+B})\textbf{u}(t) + \textbf{g}(t)
\label{eq:stiff-fode-equation}
\end{equation}
where the fractional order is chosen as $\alpha=0.5$ and final time as $T=1$, initial conditions as $u_0=(u_{01}, u_{02}, u_{03})^\top = (1,1,1)^\top$, the coefficient matrix $A$ and $B$ are:
\begin{equation}
A = \left(\begin{matrix}
    -10000 & 0 & 1\\
    -0.05 & -0.08 & -0.2\\
    1 & 0 & -1
\end{matrix}\right),\quad B = \left(\begin{matrix}
    -0.6 & 0 & 0.2\\
    -0.1 & -0.2 & 0\\
    0 & -0.5 & -0.8
\end{matrix}\right)
\end{equation}
and
\begin{equation}
g(t) = \left(\begin{matrix}
     a_1\Gamma_1t^{\sigma_1-\beta}+a_2\Gamma_2t^{\sigma_2-\beta} \\
     a_3\Gamma_3t^{\sigma_3-\beta}+a_4\Gamma_4t^{\sigma_4-\beta} \\
     a_5\Gamma_5t^{\sigma_5-\beta}+a_6\Gamma_6t^{\sigma_6-\beta}
\end{matrix}\right)-(A+B)\left(\begin{matrix}
     a_1\Gamma_1t^{\sigma_1}+a_2\Gamma_2t^{\sigma_2} + u_{01} \\
     a_3\Gamma_3t^{\sigma_3}+a_4\Gamma_4t^{\sigma_4} + u_{02} \\
     a_5\Gamma_5t^{\sigma_5}+a_6\Gamma_6t^{\sigma_6} + u_{03}
\end{matrix}\right)
\end{equation}
where $\Gamma_k=\frac{\Gamma(\sigma_k+1)}{\Gamma(\sigma_k+1-\beta)}\ (1\leq k\leq6)$, with the analytical solution of $u(t) = (a_1\Gamma_1t^{\sigma_1}+a_2\Gamma_2t^{\sigma_2} + u_{01}, a_3\Gamma_3t^{\sigma_3}+a_4\Gamma_4t^{\sigma_4} + u_{02}, a_5\Gamma_5t^{\sigma_5}+a_6\Gamma_6t^{\sigma_6} + u_{03})^\top$.

\begin{figure}[H]
    \centering
    \adjustbox{trim=4em 2em 3em 3em, clip}{
        \includegraphics[width=1.15\textwidth]{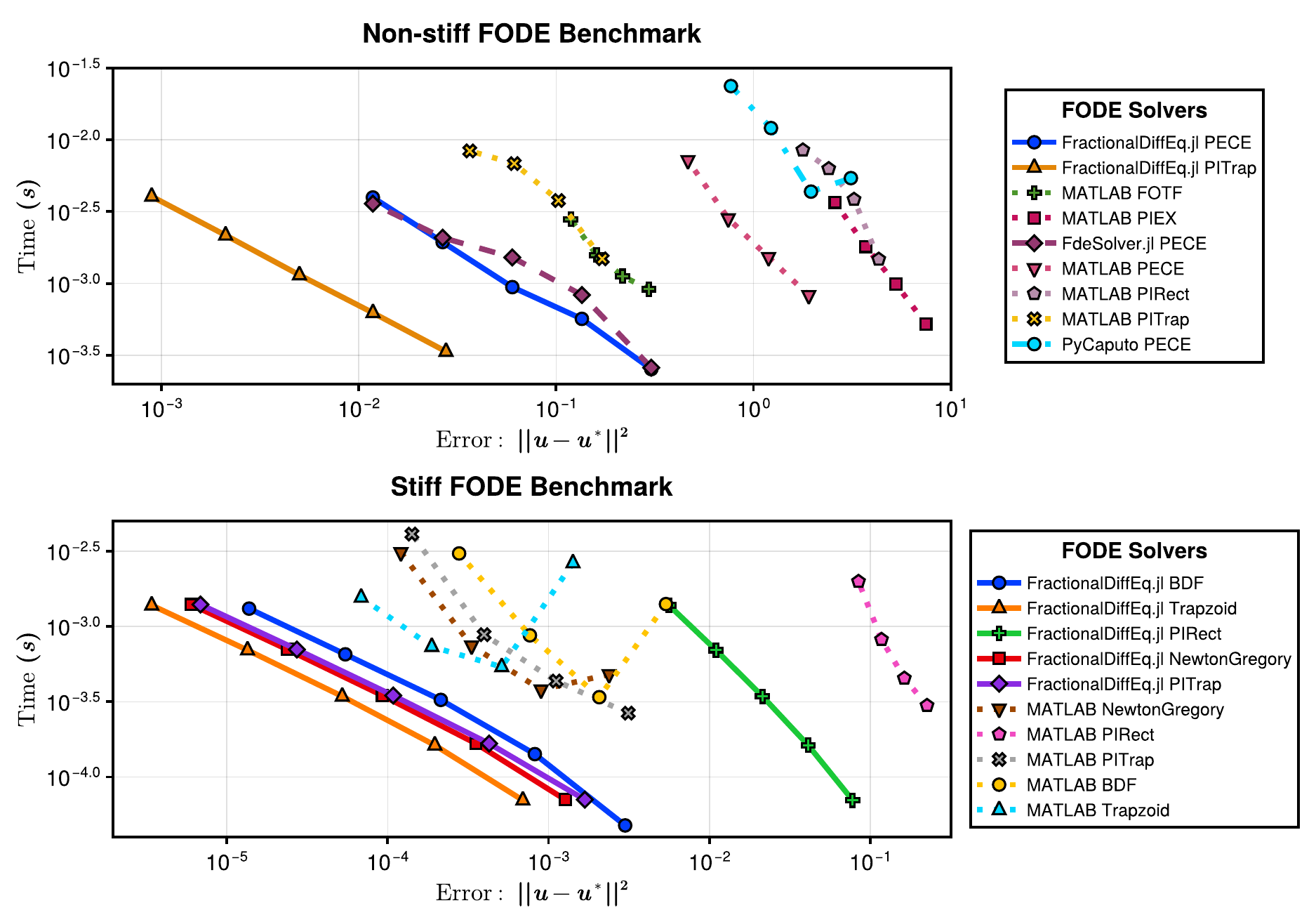}
    }
    \caption{\textbf{Work Precision Diagram for non-stiff and stiff FODE problems:}: In the benchmarking of non-stiff linear FODE problem(upper plot) from \cite{garrappa2018numerical}, solvers like \texttt{PECE}, \texttt{PITrap} from \fde~and MATLAB, \texttt{PECE} from \fdesolver~achieved convergence in considerable time, but solvers from \fde~reliably solve the nonstiff FODE problem much faster than the corresponding MATLAB and \fdesolver equivalents. On the stiff nonlinear FODE problem(lower plot) from \cite{implicit}, predictor-corrector solvers and explicit solvers failed to converge, so only the implicit methods like \texttt{BDF}, \texttt{PITrap} which are the fractional linear multistep methods from \fde~and MATLAB can successfully solve the stiff FODE problem, while all solvers from \fdesolver~and pycaputo failed to solve stiff nonlinear FODE problem at any tolerance and was removed from the plot.}\label{fig:fode-general-benchmarks}
\end{figure}

We ran the benchmarks for a step size of $h=2^{-n}, n=2,\cdots,7$. \Cref{fig:fode-general-benchmarks} shows the performance comparison of all the FODE numerical solvers from \fde, \fdesolver, pycaputo, and routines from MATLAB. Our findings indicate that in non-stiff FODE benchmarking, all solvers in \fde, \fdesolver, pycaputo, and MATLAB routines achieved convergence in considerable time, but solvers from \fde~reliably solve the non-stiff FODE problem much faster than the corresponding \fdesolver~, pycaputo and MATLAB equivalents. While in benchmarks for the stiff FODE problem from \cite{implicit}, explicit predictor-corrector solvers failed to converge, only implicit fractional linear multistep methods in \fde~and MATLAB like \texttt{BDF}, \texttt{PITrap} can successfully solve the stiff FODE problem, but all the solvers from \fdesolver~and pycaputo failed to solve the stiff nonlinear FODE problem at any tolerance.

Our empirical findings further reveal that the stiff nonlinear fractional ordinary differential equation example presented in \cite{fdesolver} is inappropriately selected, leading to erroneous conclusions regarding the solver's capabilities. And their claims that its predictor-corrector and Newton-Raphson solver can handle stiff problems is not substantiated by our experimental results. While \cite{fdesolver} claims that \fdesolver~has great reliability in various scenarios, our empirical evaluations indicate that \fdesolver~exhibits notable deficiencies in robustness and accuracy when applied to stiff fractional ordinary differential equations thereby limiting its applicability in complex modeling scenarios. Furthermore, the lack of special functions handling such as Mittag-Leffler functions in \fdesolver~and pycaputo imposes significant constraints on their ability to model fractional-order systems that inherently rely on such functions for analytical or numerical representations.

\subsubsection{Linear Multi-term FODE solvers benchmarks}\label{subsubsec:linear_multi_term_fode_benchmarks}

\begin{figure}[H]
    \centering
    \adjustbox{trim=4em 2em 3em 3em, clip}{
        \includegraphics[width=1.15\textwidth]{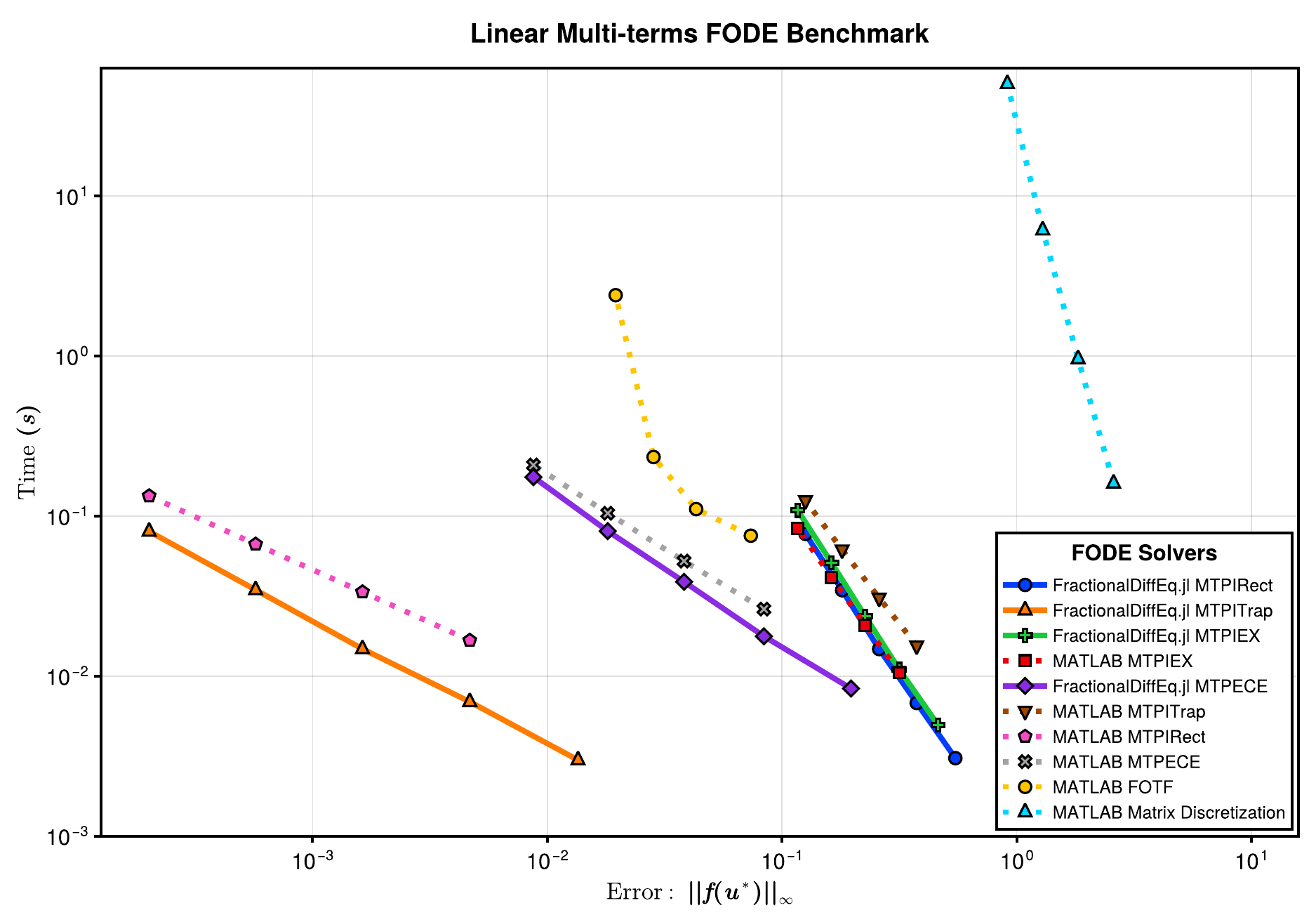}
    }
    \caption{\textbf{Work Precision Diagram of linear non-stiff multi-term fractional ordinary differential equations}: In the benchmarking of non-stiff linear problem from \cite{garrappa2018numerical}, solvers like \texttt{PITrap} and \texttt{PECE} from \fde~and MATLAB achieved good convergence in considerable time, but solvers from \fde~reliably solve the non-stiff fractional ordinary differential equations much faster than the corresponding MATLAB equivalents. There are no solvers available for multi-term FODE problems in both \fdesolver~and pycaputo.}\label{fig:linear-multiterms-benchmarks}
\end{figure}

Linear multi-term fractional ordinary differential equations (FODEs) frequently emerge in a wide range of scientific and engineering applications, including state-space representations, transfer function modeling, and anomalous diffusion processes. These equations are particularly valuable for capturing complex dynamical behaviors where classical integer-order models fail to provide accurate descriptions. To rigorously evaluate the performance of multi-term FODE solvers in \fde~, we employed the standard non-stiff benchmarking problems from \cite{XueBai+2017+1305+1312} to assess their effectiveness on non-stiff multi-term FODE problem. These benchmarking tests provide a comprehensive analysis of solver stability, accuracy, and computational efficiency under different problem conditions. The benchmarking problem from \cite{XueBai+2017+1305+1312} can be stated as:
\begin{equation}
    u'''(t)+{^C_0D^{2.5}_t}u(t)+u''(t)+4u'(t)+{^C_0D^{0.5}_t}u(t)+4u(t)=6\cos(t)
\end{equation}
with nonzero initial conditions $u(0)=1, u'(0)=1, u''(0)=-1$, and the analytical solution is $u(t)=\sqrt{2}\sin(t+\pi/4)$, with $t\in(0,100)$. We ran this benchmark for a step size of $h=2^{-n},n=2\cdots,7$. \Cref{fig:linear-multiterms-benchmarks} indicates the performance comparison between \fde~and MATLAB equivalents. Solvers from \fde~and MATLAB achieved excellent convergence on this benchmarking problem but multi-term FODE problem solvers from \fde~are faster and can achieve a higher precision.

We conducted benchmarking experiments for multi-term fractional ordinary differential equation (FODE) solvers using a step size of $ h = 2^{-n} $, where $ n = 2, \dots, 7 $. The evaluation included solvers implemented in \fde~alongside numerical routines from \cite{garrappa2018numerical}, \cite{podlubny2000matrix}, and \cite{xue2019fotf}. As illustrated in \Cref{fig:linear-multiterms-benchmarks}, the results demonstrate that the linear multi-term FODE solvers in \fde~consistently outperform their counterparts in terms of computational efficiency. Other solvers like FOTF and Matrix-Discretization methods, though they can obtain a relatively accurate results, their computational efficiency remains a significant limitation, preventing them from solving equations in a timely manner. Specifically, \fde~exhibits superior performance compared to existing methods based on FOTF, Matrix-Discretization techniques, and other MATLAB-based routines. These findings in our empirical experiments underscore the advantages of \fde~in solving complex multi-term FODEs, highlighting its potential as a robust and efficient tool for fractional-order modeling.

\section{Applications}\label{sec:applications}

In the following chapter, we delve into three exemplary application scenarios that highlighting the unique advantages of fractional ordinary differential equations in modeling complex systems. These scenarios are chosen to underscore the limitations of traditional integer-order differential equations and demonstrate how fractional calculus provides a more versatile and accurate framework for capturing the intricate dynamics of such systems.  And in these application scenarios, we utilized \fde~ for the convenient and efficient simulation, further demonstrating the capability and applicability of \fde~.

The first scenario in \Cref{subsec:fermentation} involves a state-space model with inherent memory effects or hereditary properties, where the current state depends not only on the immediate past but also on the entire history of the system. Such behavior is common in many physical, biological, and engineering contexts, making FODEs particularly well-suited for describing these dynamics. By incorporating fractional derivatives, which naturally account for non-local interactions over time, these models achieve a level of precision and realism unattainable with integer-order approaches.

The second scenario in \Cref{subsec:harmonic-oscillator} entails a harmonic oscillator with fractional orders, represented as a linear multi-term fractional ordinary differential equation. \Cref{subsec:harmonic-oscillator} explores how the behavior of the oscillator evolves in fractional order as the oscillator system varies, providing insights into the impact of the fractional dynamics on system stability and response characteristics.

Through these examples, we aim to illustrate the practical relevance and superior modeling capabilities of fractional differential equations. By leveraging \fde, these models open new avenues for understanding and simulating complex phenomena that would otherwise remain inadequately described. This chapter serves as a testament to the growing importance of fractional calculus in addressing the challenges of modern science and engineering, and how \fde~plays an important role in the real-world simulation of fractional differential equations.

\subsection{Fractional-order fermentation models}\label{subsec:fermentation}

Recent studies have increasingly highlighted the growing significance of using fractional differential equations in modeling real-world phenomena, particularly in applications such as bio-kinetic analysis and epidemiology analysis \cite{da2024fractional, barbero2024brief}. The adoption of fractional-order models is driven by their ability to more accurately represent systems exhibiting memory effects and anomalous diffusion, which are often encountered in biological, ecological, and epidemiological processes.

The fractional fermentation model proposed by \cite{toledo2014fractional} is an analogy of the classical Lotka-Volterra model, incorporating three key state variables: the biomass ($B$), the product ($P$), and the substrate ($S$). The fractional-order tequila fermentation model \Cref{eq:lotka-volterra} assumes that the biomass growth is driven solely by the available of the substrate, with the rate of biomass accumulation being proportional to the interaction between biomass and substrate($BS$ in \Cref{eq:lotka-volterra}). This interaction reflects the fundamental principle that microbial growth is sustained by substrate consumption, aligning with the classical Monod kinetics \cite{kovarova1998growth} but with an added fractional-order dynamic to capture complex temporal dependencies.

Let's consider the fermentation process of tequila production in \cite{arellano2007unstructured} and \cite{toledo2014fractional}, 
\begin{equation}
    \begin{aligned}
        &{_0^C}D_t^{\alpha_1}B=k_cBS-k_mB\\
        &{_0^C}D_t^{\alpha_2}S=-k_sBS\\
        &{_0^C}D_t^{\alpha_3}P=k_pBS
    \end{aligned}
    \label{eq:lotka-volterra}
\end{equation}
where $B$ is the biomass concentration, $S$ is the substrate concentration, and $P$ is the product concentration(here is the ethanol), the kinetics parameters and their corresponding units in \Cref{eq:lotka-volterra} are described in \Cref{tab:quantities_in_bioethanol}.
\begin{table}[H]
    \centering
    \begin{tabular}{ccccc}
    \toprule
         Notation & Description & Value & Units \\
         \midrule
         $k_m$ & mortality rate & 0 & $1/s^\beta$  \\
         $k_c$ & grow rate of the biomass & 0 & $1/(g\cdot s^\beta)$  \\
         $k_k$ & consumption rate of the substrate & 0 & $1/(g\cdot s^\beta)$  \\
         $k_p$ & formation rate of the final product ethanol & 0 & $1/(g\cdot s^\beta)$  \\
         \bottomrule
    \end{tabular}
    \caption{Description of the parameters used in model \Cref{eq:lotka-volterra} and their corresponding units.}
    \label{tab:quantities_in_bioethanol}
\end{table}

While the ordinary differential equations model for the tequila fermentation process is given in \cite{arellano2007unstructured, toledo2014fractional}:
\begin{equation}
\begin{cases}
    &\frac{dx_1}{dt} = (\mu-D)x_1\\
    &\frac{dx_2}{dt} = -(Y_{xs} - Y_{ps})\mu x_1-m_sx_1 - (-S_0+x_2)D\\
    &\frac{dx_3}{dt} = \alpha\mu x_1 - Dx_3
\end{cases}
\label{eq:fractional_tequila}
\end{equation}
where $Y_{xs}$ and $Y_{ps}$ are the dimensionless biomass and product yield coefficients respectively, $S_0$ is the feed substrate concentration, $m_s$ is the maintenance coefficient, and $\mu$ is given by $\mu=\frac{\mu_{max}x_2}{k_s+x_2+x_2^2k_1}(1-x_3k_p)$.
\begin{table}[H]
    \centering
    \begin{tabular}{c|c|c}
    \toprule
         \diagbox[width=2.5cm, height=1.3cm]{Estimations}{Models} & Integer-order model & Fractional-order model \\
         \midrule
         $k_c$ & 0.003027 & $0.004265$  \\
         $k_m$ & 0.000499 & $0.000499$  \\
         $k_s$ & 0.050000 & $0.055166$  \\
         $k_p$ & 0.015999 & $0.015805$  \\
         $\alpha_1$ & - & $0.775347$    \\
         $\alpha_2$ & - & $0.873674$    \\
         $\alpha_3$ & - & $0.976698$    \\
         \hline
         RMSD & 0.036507 & 0.012339 \\
         \bottomrule
    \end{tabular}
    \caption{Estimated optimized parameters and fractional orders used in integer-order model and fractional-order model \Cref{eq:lotka-volterra}. The quantitative assessment of the model fidelity was extended through rigorous computation of root mean square deviation metrics of the different model to systematically quantify absolute errors.}
    \label{tab:estimation_in_bioethanol}
\end{table}

\begin{figure}[H]
        \centering
        \includegraphics[width=1.0\textwidth]{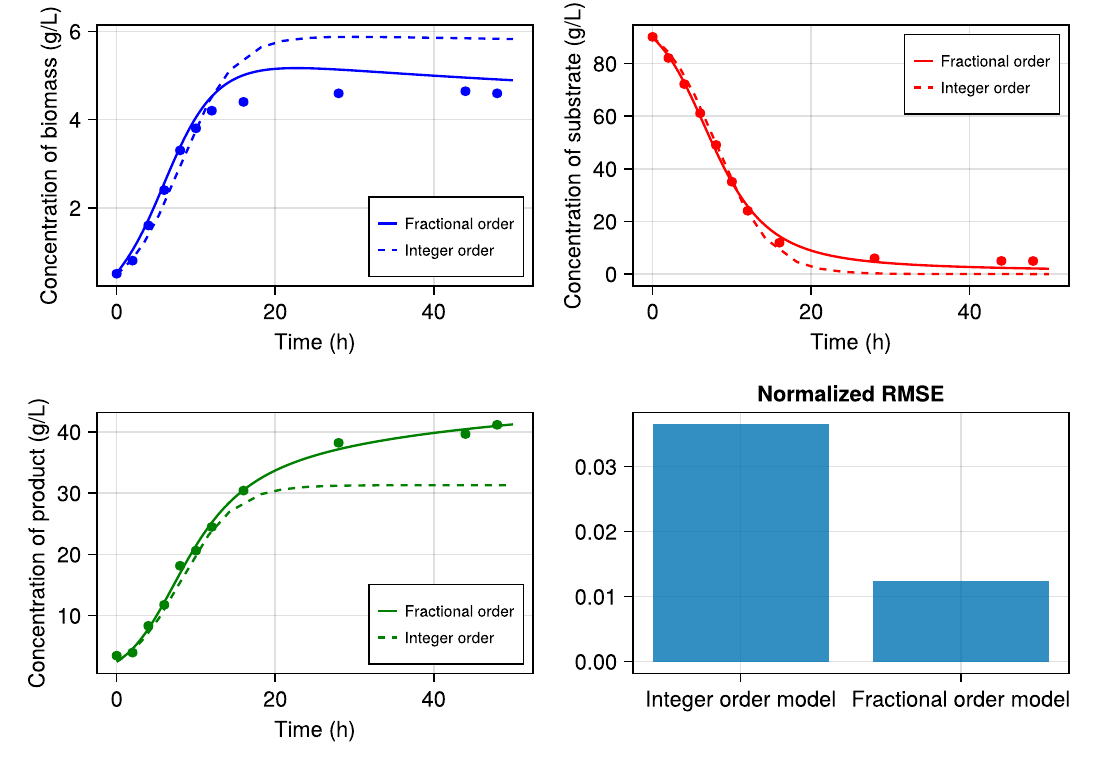}
    \caption{The comparison between real data from experiments as retrieved from \cite{toledo2014fractional} and the estimation from the fractional-order tequila fermentation model \Cref{eq:fractional_tequila} with integer-order derivatives(Integer-order model shown by dash line) and incommensurate fractional-order derivatives(Fractional-order model shown by solid line). The circle indicates the real data during a time span of 48h with sampling performed every 2 h during the first 12 h of fermentation, then every 4 h during the following 36 h. The specified error are based on the root mean square error (RMSE) where fewer value represents less residual variance and better fitting results. The final fitted quantities are listed in \Cref{tab:quantities_in_bioethanol}}\label{fig:fermentation}
\end{figure}

The kinetic parameters of the fractional-order and integer-order tequila fermentation model were both estimated from nonlinear curve fitting with Optim.jl \cite{mogensen2018optim}. This approach ensures precise parameter estimation by minimizing the error between the analytical model predictions in \Cref{eq:fractional_tequila} and \Cref{eq:fractional_tequila} with experimental data, allowing for a more accurate characterization of the fermentation dynamics. The experimental data are collected from experiments in \cite{toledo2014fractional} and \cite{arellano2007unstructured} where sampling was performed every 2 h during the first 12 h of fermentation, then every 4 h during the following 36 h. We utilized the LBFGS optimization solver which is a quasi-Newton method from \cite{nocedal1999numerical} to minimize the residual during the nonlinear fitting of parameters for model \Cref{eq:lotka-volterra} and model \Cref{eq:fractional_tequila}, find the optimal value to best fit the integer-order and fractional-order model to the experimental data. All the fitted values are listed in \Cref{tab:estimation_in_bioethanol}. The root mean square errors in \Cref{fig:fermentation} is calculated using StatsBase.jl \cite{dahua_lin_2024_14438491} defined as RMSE=$\sqrt{\displaystyle\sum_{i=1}^n\frac{(\hat{y_i}-y_i)^2}{n}}$. We use the computed root mean square error to evaluate the fitting performance of different model.

While Caputo fractional derivatives are defined by integrals, they are basically non-local operators and are able to describe dynamical models whose dynamical process involves memory effect. From the final results of experiments in \Cref{fig:fermentation}, the biological process of fermentation for tequila production possesses a physicochemical nature that involves the history of bio-reactions, by utilizing the fractional-order fermentation model, we can better describe the underlying mathematical model compared to using integer-order model with a small number of parameters and a much simpler numerical expression, and the predictive ability if more accurate than ordinary models. The parameter estimation from \Cref{fig:fermentation} indicates that fractional-order models exhibit superior capability in capturing the system dynamics with relatively small error compared to the integer-order counterparts, especially in simulations like the fermentation model which involves memory effects. Such memory effects indicate that past products and substrates can play an important role in the current reaction, and their contributions to the whole fermentation process are negligible.

\subsection{Fractional harmonic oscillator}\label{subsec:harmonic-oscillator}
Fractional multi-term FODE provide a new way to model complex dynamical systems characterized by memory effects and anomalous diffusion processes. \fde~offer an efficient and robust computational framework for solving linear multi-term fractional differential equations, addressing key challenges such as numerical stability, accuracy, and computational efficiency in dealing with linear multi-term FODEs. By leveraging advanced discretization techniques, optimized iterative solvers, and efficient memory management strategies, our approach significantly enhances the practical applicability of multi-term FODEs in large-scale simulations and real-world modeling scenarios.

Consider the fractional harmonic oscillator system which is actually a linear multi-term fractional-order initial value problem:
\begin{equation}
\begin{aligned}
    &u''+\alpha u'+\beta{^CD^\theta}u+\gamma u=0 \\
    &u(0)=-0.5,\ u'(0)=-0.5
    \end{aligned}
    \label{eq:fractional-harmonic-oscillator}
\end{equation}
where $\alpha$ is the viscous damping coefficient and $\beta$ is the friction coefficient.
more background
The fractional harmonic oscillator has a characteristic equation of:
\begin{equation}
    \Delta(s;\alpha,\beta,\gamma,\theta) = s^2+\alpha s+\beta s^\theta+\gamma=0
\end{equation}
hence the oscillator system is asymptotically stable when $\alpha>0,\beta>0$ and $\gamma>0$, as established by the stability criterion in \cite{brandibur2021stability}. This ensure that all eigenvalues of the system lie within the stability region, leading to the decay oscillation over time. However, when either the coefficient $\alpha$ or $\beta$ become negative, the stability behavior of the system changes significantly. In this case, solving the characteristic equation \Cref{eq:fractional-harmonic-oscillator} reveals that the stability of the oscillator is no longer solely determined by the sign of the coefficients but also by the fractional order $\theta$, which has a critical value of $\theta^\star=0.81695$ determining the stability of the oscillator.

\begin{figure}[H]
        \centering
        \includegraphics[width=0.6\textwidth]{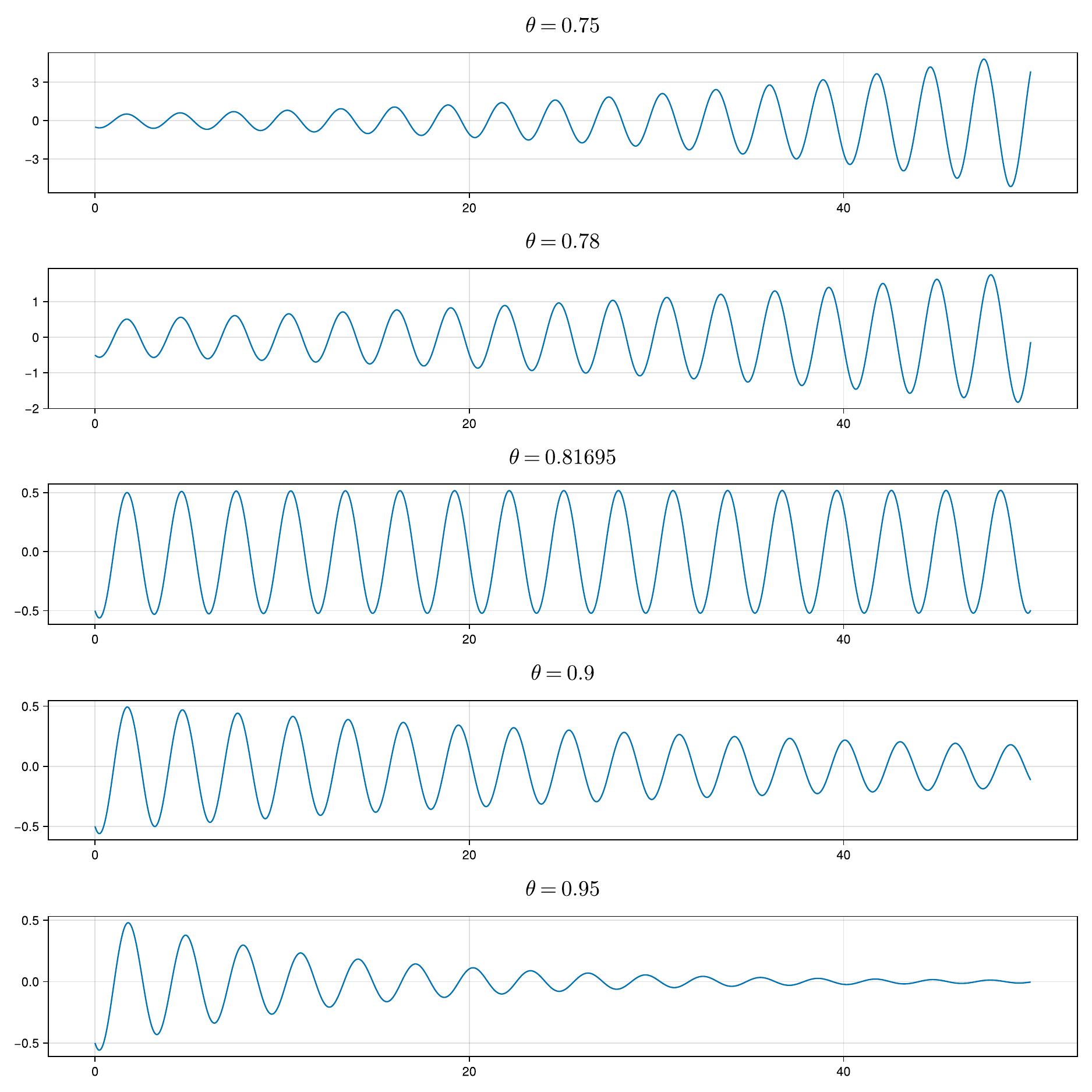}
    \caption{Different behavior of fractional harmonic oscillator in different fractional order $\theta$ when the parameters $\alpha, \beta, \gamma$ are chosen as $\alpha=-1,\beta=1.2,\gamma=4$. When $\theta$ is $0.81695$, which is the critical point of \Cref{eq:fractional-harmonic-oscillator}, the oscillator would maintain a periodic oscillation.}\label{fig:multiterms-application}
    \vspace{-2em}
\end{figure}
The \Cref{fig:multiterms-application} illustrates the dynamical behavior of the fractional oscillator system \Cref{eq:fractional-harmonic-oscillator}. Specifically, when the fractional order $\theta<\theta^\star$ the oscillator system exhibits instability,characterized by an exponential growth in the amplitude of oscillations over time, ultimately resulting in a divergent response. This phenomenon suggests that for subcritical fractional orders, the system lacks sufficient damping mechanisms to counteract perturbations, leading to unbounded energy accumulation. Conversely, when $\theta \geq \theta^\star$, the system transitions into a stable regime, where oscillations either asymptotically decay to a trivial equilibrium state or persist in a bounded periodic manner. This stability threshold $\theta^\star$ thus plays a critical role in governing the long-term behavior of the system, delineating the boundary between instability-induced divergence and controlled, physically meaningful oscillations. The observed bifurcation behavior underscores the sensitivity of fractional-order dynamical systems to parameter variations, emphasizing the necessity of precise fractional order selection in practical applications.

\section{Discussion}
\label{sec:discussion}

In this work, we present a high performance package for numerically solve fractional differential equations within the Julia language and SciML ecosystem. In \Cref{subsec:benchmarks}, we demonstrated \fde~superior performance comparing counterparts, especially in \Cref{subsubsec:fode_benchmarks} benchmarks on stiff fractional ordinary differential equations, \fde~ outperforms \fdesolver~ and other solver routines available from MATLAB and Python, proves its superiority in solving difficult fractional ordinary differential equations.

With the utilization of Julia's multiple dispatch and appropriate type stable code, we observed a huge speedup compared to the MATLAB routines and another implementation in \fdesolver~and Python pycaputo package, further differentiating \fde~ in complex fractional-order modeling area for its great capabilities and high performance. \fde~share the same interface and functionalities with SciML differential equations solvers ecosystem, which means both the fundamental utilities and advanced analyzing tools from SciML can be applied to \fde~as well, broaden the coverage of SciML differential equations solving suite and provide \fde~further modeling development such as acausal modeling with ModelingToolkit.jl \cite{ma2021modelingtoolkit}, sensitivity analysis with SciMLSensitivity.jl \cite{rackauckas2020universal}, surrogate modeling with Surrogates.jl \cite{ludovicobessi_2024_12571719} and neural fractional differential equations with DiffEqFlux.jl \cite{rackauckas2019diffeqflux} and NeuralPDE.jl \cite{zubov2021neuralpde}.

Through comprehensive and exhaustive benchmarking, we conducted an in-depth comparison of the computational performance of solving fractional ordinary differential equations using \fde~against MATLAB, Python, and its Julia-based alternatives. The results demonstrated a remarkable improvement in performance when employing \fde~ and highlights the superior efficiency and optimization of \fde~in handling complex fractional differential equations problem. In particular, the enhanced performance stems from its advanced algorithms, which are specifically designed to tackle the computational challenges associated with fractional calculus. The benchmarks revealed that \fde~not only surpasses MATLAB in execution speed but also outperforms existing Julia and Python solutions, positioning it as a leading tool for researchers and practitioners working on fractional differential equations. This breakthrough underscores the potential of \fde~to streamline computational workflows and enable faster and more accurate simulations in various scientific and engineering applications.

Additionally, in contrast to the limited scope of fractional differential equations types other solver packages are able to cover, \fde~encompasses significantly more types of fractional differential equations including multi-term fractional ordinary differential equations, while in other alternatives such as \fdesolver~and pycaputo, there are no such built-in functionalities to address these complex multi-term fractional systems, thereby restricting their applicability in more advanced modeling scenarios. By providing robust fractional ordinary differential equations solvers, \fde~can bridge the gap between research in numerical solvers for fractional differential equations and real-world numerical software, illustrating its significance in computational research and application development.

The applications of \fde~in various kinds of fractional differential equations in \Cref{sec:applications} further differentiate its distinctive capability in efficiently solving complex fractional differential equations across diverse problem domains. These applications demonstrate the robustness, computational efficiency, and adaptability of \fde~in handling both linear and nonlinear fractional systems, as well as multi-term and stiff fractional differential equations. By leveraging advanced numerical schemes and optimized implementations, \fde~exhibits superior performance in terms of accuracy, convergence rates, and computational scalability compared to conventional solvers. The demonstrated versatility of \fde~across different application scenarios underscores its potential as a powerful computational tool for solving fractional-order models arising in science and engineering.

\section{Conclusion}
\label{sec:conclusion}

Solving fractional differential equations is becoming a fundamental challenge that arises across various scientific domains. This paper presented \fde, a high-performance and robust open-source solver for fractional differential equations implemented natively in the Julia programming language. To provide easy and convenient user interfaces, \fde~follow the design pattern of SciML DifferentialEquations.jl \cite{rackauckas2017differentialequations} which deployed a "define problem-solve problem" semantics to handle ODE, SDE, BVP, and DAEs easily. To help users with a better understanding of the usage of the \fde~package, thorough documentation and numerical examples are provided in the accompanying user manual and website, further facilitating researchers' deployment of fractional-order numerical models and being indispensable tools for numerical FDE solving. With thorough usage examples, exhaustive numerical experiments on benchmark problems, and real-world applications, we have demonstrated the superior capabilities of \fde~compared to other existing software tools.

Key strengths of \fde~include its flexible unified API for rapidly experimenting with different solver options, and exceptional efficiency built on top of the Julia language which is a high-level dynamical language utilizing the multiple dispatch paradigm for high-performance numerical computing. \fde~provides mature numerical solvers for more types of fractional differential equations, including system of fractional ordinary differential equations and linear multi-term fractional ordinary differential equations. These features enable \fde~to reliably solve challenging fractional differential equations, including cases where standard solvers fail while attaining high performance.

Future development for \fde~include integrating NonlinearSolve.jl \cite{pal2024nonlinearsolve} for robust nonlinear solvers in fractional linear multistep methods. Massively parallel solving differential equations across GPU threads could solve fractional differential equations in parallel across GPU threads \cite{utkarsh2024automated, besard2018effective}, hence accelerating large complex fractional differential equations solving with \fde~is a big step towards large-scale simulation of complex fractional ordinary differential equations.

Ultimately, \fde~not only advances fractional-order modeling and simulations by providing multiple high-performance and robust numerical solvers for different types of fractional ordinary differential equations, but also provides simple and unifying software interface, making itself a dispensable and powerful tool for researchers and practitioners in various scientific and engineering fields. \fde~ helps bridge the gap between theoretical research and real-world application on fractional calculus and fractional differential equations by providing a robust and scalable solvers suite for simulating complex systems. Its intuitive API and high-performance algorithms enable researchers to easily model, analyze, and solve fractional differential equations, thus accelerating the translation of theoretical insights into practical solutions. Its versatility and efficiency enable the exploration of complex dynamical systems that incorporate memory effects, thereby broadening the scope of applications in areas such as control systems, bioengineering, finance, and material science.

\section*{Acknowledgment}
\label{sec:acknowledgement}

The author would like to sincerely acknowledge Zhejiang University and Shandong University for providing invaluable academic resources and a stimulating research environment that greatly facilitated this study. We are also deeply grateful our gratitude to the anonymous reviewers for their insightful comments and constructive feedback, which helped improve the quality of this manuscript.


\bibliographystyle{unsrtnat}
\bibliography{references}

\end{document}